\documentclass{amsart}
\usepackage{setspace, amsmath, amsthm, amssymb, amsfonts, amscd, epic, graphicx, ulem, dsfont}
\usepackage[T1]{fontenc}
\usepackage{multirow}
\usepackage{bbm}
\usepackage{enumerate}

\makeatletter \@namedef{subjclassname@2010}{
  \textup{2020} Mathematics Subject Classification}
\makeatother

\newtheorem{thm}{Theorem}[section]

\theoremstyle{remark}

\theoremstyle{definition}

\newtheorem{exa}[thm]{\textbf{Example}}

\newcommand{\R}{\mathbb{R}}

\begin{document}

\title{Simple examples of non closable paranormal operators}
\author[M. H. MORTAD]{Mohammed Hichem Mortad}

\thanks{}
\date{}
\keywords{Closable operator. Paranormal operator. Hilbert space.
Fourier transform}

\subjclass[2010]{Primary 47B20. Secondary 47A05. }

\address{Department of
Mathematics, University of Oran 1, Ahmed Ben Bella, B.P. 1524, El
Menouar, Oran 31000, Algeria.\newline {\bf Mailing address}:
\newline Pr Mohammed Hichem Mortad \newline BP 7085 Seddikia Oran
\newline 31013 \newline Algeria}

\email{mhmortad@gmail.com, mortad.hichem@univ-oran1.dz.}

\begin{abstract}In this note, we give an example of a densely defined non-closable
paranormal operator. Then, we give another example of a densely
defined closable paranormal operator whose closure fails to be
paranormal. These two examples are simpler than those which first
appeared in \cite{Daniluk-paranormals-non-closable}. It is worth
noticing that our first example tells us that the adjoint of a
densely defined paranormal operator may have a trivial domain.
\end{abstract}

\maketitle

\section{Introduction}

First, we assume readers familiar with basic notions and results
about linear unbounded operators, as well as matrices of non
necessarily bounded operators. Two useful references are
\cite{SCHMUDG-book-2012} and \cite{tretetr-book-BLOCK} respectively.
Some basic knowledge of the $L^2(\R)$-Fourier transform (denoted by
$\mathcal{F}$ throughout the paper) is also needed.

Recall that a linear operator $B$ with a domain $D(B)$ is called an
\textbf{extension} of another linear operator $A$ (with a domain
$D(A)$), and we write $A\subset B$, if
\[D(A)\subset D(B) \text{ and } \forall x\in D(A):~Ax=Bx.\]

A linear operator $A$ is said to be densely defined if
$\overline{D(A)}=H$. Say that a linear operator $A$ is closed if its
graph is closed in $H\oplus H$. A linear operator $A$ is called
closable if it has a closed extension, the smallest (w.r.t.
"$\subset$") closed extension is called its closure and is denoted
by $\overline{A}$. It is also known that a densely defined linear
operator $A$ is closable if and only if $D(A^*)$ is dense in $H$.

Recall also that a linear operator $A:D(A)\subset H\to H$ (where $H$
is a Hilbert space) is said to be \textbf{paranormal} if
\[\|Ax\|^2\leq \left\|A^2x\right\|\|x\|\]
for all $x\in D(A^2)$. This is clearly equivalent to $\|Ax\|^2\leq
\|A^2x\|$ for all \textit{unit} vectors $x\in D(A^2)$.

Recall that $C_0^{\infty}(\R)$ denotes here the space of infinitely
differentiable functions with compact support. The following result,
whose proof relies upon the Paley-Wiener's theorem, is well known.

\begin{thm}\label{Fourier trans f, hat f C_0inf f=0 thm}
If $f\in C_0^{\infty}(\R)$ is such that its Fourier transform $\hat
f\in C_0^{\infty}(\R)$, then $f=0$.
\end{thm}

Note in the end that we tolerate the abuse of notation $A(x_1,x_2)$
instead of $A\left(
     \begin{array}{c}
       x_1 \\
       x_2 \\
     \end{array}
   \right)$ in the case of a matrix of operators $A$.

(Non necessarily bounded) Quasinormal, subnormal and hyponormal
operators are all closable. The closability of hyponormal operators
(which is obvious) suffices for the closability of other two classes
thanks to the usual inclusion
\[\text{Quasinormal}\subset\text{Subnormal}\subset \text{Hyponormal}.\]

However, as is known, every hyponormal operator is paranormal and so
it is natural to ask whether paranormal operators are closable?
Another natural question is whether a closable paranormal operator
must have a paranormal closure? These two questions were first
answered negatively by A. Daniluk in
\cite{Daniluk-paranormals-non-closable}, after having remained open
for some while. Daniluk's counterexamples were not too complicated
but they did require quite a few results before getting to the final
counterexamples. The main aim of this paper is to provide simpler
counterexamples than Daniluk's.

\section{Main Counterexamples}

First, we give an example of a densely defined linear operator $T$
with domain $D(T)$ such that $T$ is a non closable paranormal
operator. Incidentally, it may well happen that a densely defined
paranormal operator $T$ satisfies $D(T^*)=\{0\}$ (unlike other
classes such as symmetric and quasinormal operators).

\begin{exa}
Let $T$ be a densely defined operator $T$ such that
\[D(T^2)=D(T^*)=\{0\}.\]
Such example appeared in Proposition 2.2 in \cite{Mortad-TRIVIALITY
POWERS DOMAINS}. For readers convenience recall that there
$T:=A^{-1}B$ where $A$ and $B$ are defined by
\[Af(x)=e^{\frac{x^2}{2}}f(x)\]
on $D(A)=\{f\in L^2(\R):~e^{\frac{x^2}{2}}f\in L^2(\R)\}$ and
$B:=\mathcal{F}^*A\mathcal{F}$ (see \cite{Mortad-TRIVIALITY POWERS
DOMAINS} for all the details).

Then for $x\in D(T^2)=\{0\}$
\[\|Tx\|^2=\|T^2x\|\|x\|=0,\]
i.e. $T$ is trivially paranormal. Since also $D(T^*)=\{0\}$, $T$
cannot be closable.
\end{exa}

Next, we exhibit a closable densely defined paranormal operator $T$
such that its closure $\overline{T}$ is not paranormal.

\begin{exa}
Let $\mathcal{F}$ be the usual $L^2(\R)$-Fourier transform and let
$A$ be the restriction of $\mathcal{F}$ to the dense subspace
$C_0^{\infty}(\R)$. Now, define
\[T=\left(
      \begin{array}{cc}
        0 & A \\
        \frac{1}{2}A & 0 \\
      \end{array}
    \right)
\]
on $D(T)=C_0^{\infty}(\R)\oplus C_0^{\infty}(\R)$. Since $D(T)$ is
clearly dense in $L^2(\R)\oplus L^2(\R)$, $T$ is densely defined.
Since
\[\left(
      \begin{array}{cc}
        0 & A \\
        \frac{1}{2}A & 0 \\
      \end{array}
    \right)\subset \left(
      \begin{array}{cc}
        0 & \mathcal{F} \\
        \frac{1}{2}\mathcal{F} & 0 \\
      \end{array}
    \right)\]
and the latter is clearly closed (in fact everywhere defined and
bounded), $T$ is closable.

Since
\[D(A^2)=\{f\in C_0^{\infty}(\R):\mathcal{F}f\in C_0^{\infty}(\R)\}=\{0\},\]
it results that $D(T^2)=\{(0,0)\}$. So $T$ is paranormal (as in the
above argument). It only remains to show that $\overline{T}$ is not
paranormal. Clearly,
\[\overline{T}=\left(
      \begin{array}{cc}
        0 & \mathcal{F} \\
        \frac{1}{2}\mathcal{F} & 0 \\
      \end{array}
    \right)\]
(remember that $\mathcal{F}$ is the everywhere defined
$L^2(\R)$-Fourier transform). Hence
\[\overline{T}^2=\left(
      \begin{array}{cc}
        0 & \mathcal{F} \\
        \frac{1}{2}\mathcal{F} & 0 \\
      \end{array}
    \right)\left(
      \begin{array}{cc}
        0 & \mathcal{F} \\
        \frac{1}{2}\mathcal{F} & 0 \\
      \end{array}
    \right)=\left(
      \begin{array}{cc}
        \frac{1}{2}\mathcal{F}^2& 0 \\
        0& \frac{1}{2}\mathcal{F}^2 \\
      \end{array}
    \right)\]
(defined now on $L^2(\R)\oplus L^2(\R)$). If $\overline{T}$ were
paranormal, it would ensue that
\[\forall (f,g)\in L^2(\R)\oplus L^2(\R):\|\overline{T}(f,g)\|^2\leq \|\overline{T}^2(f,g)\|\|(f,g)\|.\]
In particular, this would be true for some $(0,g)$ with
$\|g\|_2\neq0$, where $\|\cdot\|_2$ denotes the usual
$L^2(\R)$-norm. That is, we would have
\[\|\mathcal{F}g\|_2^2\leq \frac{1}{2}\|\mathcal{F}^2g\|_2\|g\|_2.\]

By invoking the Plancherel's theorem $\|\mathcal{F}g\|_2=\|g\|_2$,
and by the general theory, $\mathcal{F}^2g(x)=g(-x)$. Therefore, the
previous inequality would become
\[\|g\|_2^2\leq \frac{1}{2}\|g\|_2\|g\|_2\text{ or merely } 1\leq \frac{1}{2}\]
which is absurd. Accordingly, $\overline{T}$ is not paranormal, as
wished.
\end{exa}

\end{document}